# The Ethical Turn in Mathematics Education

## Dennis Müller[1]

Institute of Mathematics Education

University of Cologne, Germany

**Abstract:**

This article analyzes the emerging ethical turn in mathematics education, arguing that it is a nuanced extension of the sociopolitical turn. While sociopolitical studies of mathematics have highlighted systemic issues and group concerns (e.g., equity, diversity, exclusion), the newer scholarship on ethics in mathematics presents a sharpened focus on the individual responsibility of learners, teachers, and mathematicians by explicitly engaging with philosophical ethics. We analyze key themes of the discourse, including the tension between "doing good" and "preventing harm," and present various philosophical foundations from which scholars have engaged with ethics: Levinas, non-Western perspectives, and pragmatism. We show that the ethical turn holds significant implications for training teachers, including self-reflection, responsibility towards the Other, historical and philosophical awareness, the role of mathematics in society, individual flexibility, cultural sensitivity, and courage to navigate the complex reality of today's mathematics classrooms. The article is designed to also serve as an introduction to ethics in mathematics education.

---

[1] For feedback or inquiries please send an email to dm782@cantab.ac.uk.

# Table of Contents

# 1. <u>Introduction</u>

As mathematics education has transitioned from a discipline largely concerned with teaching and learning mathematics to researching topics in a broader social and political context (Gutiérrez, 2013), the range of critical mathematics scholarship has widened accordingly in recent decades. Building on wide arrays of postmodern[2] and poststructuralist[3] scholarship, critical educators foreground their writings in the interwovenness of knowledge, power, identity, and social discourses (Müller, 2024), aiming not only to uncover but also to overcome the structures that privilege a select group of individuals (e.g., transformative research, as outlined by Skovsmose & Borba, 2004), empower and not just study disadvantaged groups of students (e.g., Gutiérrez, 2008), but foster participation in society (Aguilar & Zavaleta, 2012, p. 6). In short, "mathematics teaching is [seen as] democratic education," as Hannaford (1998) put it. Gutiérrez (2013, pp. 51–55) identifies five key benefits of considering mathematics education from a sociopolitical perspective. These benefits include transcending essentialized group identities, moving beyond narratives of victimization towards empowerment, redefining qualitative mathematics education, challenging hierarchies by recognizing the complexity of the learners' and teachers' identities, and deconstructing traditional (modern) views of mathematics.

---

[2] Our personal view of Postmodernism largely follows Walshaw 's (2004) fairly general perspective of "asking questions about modernist thinking and acting, and propose conceptual designs to investigate an increasingly complex, plural and uncertain world"(ibid., pp. 1-2).

[3] Similarly, our understanding of Poststructuralism can be compared to Harcourt 's (2007) description as "critical reasoning that focuses on the moment of slippage in our systems of meaning […] [when] ethical choices are made," (ibid., p. 1). We particularly follow the idea that "strongly held rational beliefs in certain theories or premises [or structures] rest on a leap of faith" (ibid., p. 21). For us, particular importance lies the latent tension between the inherently localised ethics and the believe in universal, neutral and pure of mathematics, which we understand to be explored implicitly and explicitly by many of the authors reviewed in this article.

Now Ernest (2024a) has suggested the emergence of further developments, claiming that while the "widespread 'ethical turn' in philosophy, the humanities and social sciences [had long] left the research field of mathematics behind"(ibid., p. 3), "it is […] now reaching mathematics and mathematics education"(ibid., p. 4). But Ernest was not the first to notice or advocate for this turn. In the context of Australian mathematics education, Vale et al. (2016) have previously described a similar shift from psychology to sociopolitics to ethics. There, many scholars naturally did not give up their past engagement with critical pedagogy but rather used ethics as another tool in their repertoire to analyze the problems in today's mathematics classroom. We will see that this is also true here. For example, Register et al. (2021) see a need to consider ethics in mathematics education, as the current studies of social justice are often only focused on underprivileged groups. Moreover, while Müller (2024) draws a distinction between ethics in mathematics and mathematics for social justice, he sees sufficiently many similarities between them to be naturally connected and co-taught.

As sociopolitical and ethical issues are deeply entangled and mutually informative in mathematics education, our reading of the ethical turn is not a "turning away" from sociopolitical issues but a further nuance. Indeed, many foundational works addressing sociopolitical issues like equity, power, and inclusion ask fundamental ethical questions about the role of mathematics and education in society, even if they do not explicitly frame all their arguments through philosophical ethics (e.g., the ZDM special issue on *mathematics, peace and ethics* (D'Ambrosio & Marmé, 1998)). On the other hand, the common ethical questions asked by authors may also be rooted in social justice concerns, and sociopolitical questions can be analyzed from an explicitly ethical perspective (Chen, 2024). As a concrete example, the ethics of assessment, with a particular focus on marginalization, is studied by Bagger and Nieminen (2022), including through the lens of ethical dilemmas and professionalism (Bagger, 2024). Other examples include the tensions between equity and quality when explored from the perspective of ethics (Atweh, 2011).

We will show that this latest turn seeks to engage more explicitly with philosophical ethics, often extending the work initiated by the sociopolitical scholarship. To analyze the ethical turn, we consider both mathematics education at the school and university level. We believe that making a distinction here and only considering one or the other would restrict the analysis unnecessarily. Many scholars write about both issues and regularly cross the boundary in between. Additionally, with its immense authority (Ernest, 2024b), mathematical practices also transcend these boundaries, creating new and destroying old orders (cf. Ali, 2024).

The interconnectedness of different discourses becomes clear when we look at the connectivity graph of the bibliography of this article (Figure 1). In this graph, we not only see how arguments and scholars transcend the various boundaries but also how central fundamental works on the sociopolitical nature of mathematics education are to the ethical turn. We can also observe an increasing number of publications on the topic in recent years, with an increase in interconnectivity. This preemptively suggests that the ethical turn may not just be a turn for individual authors but is indeed beginning to show signs of a discursive life of its own. Additionally, Figure 2 foreshadows a timeline of the various foci within the ethics in mathematics education literature.

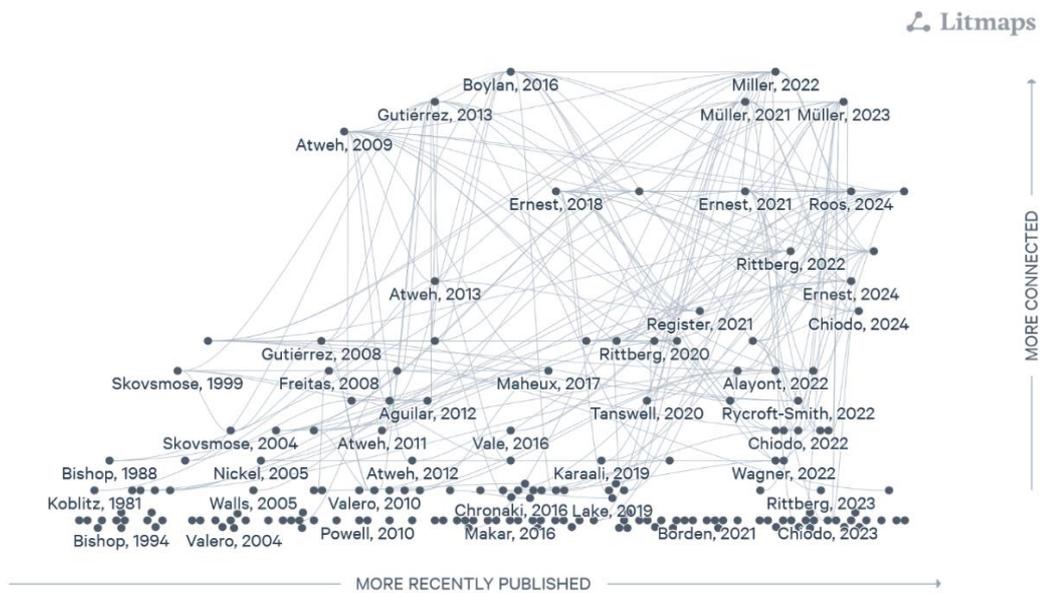

Figure 1: Connectivity graph of the bibliography[4]

| Pre-2000 | Limited explicit focus on ethics |
| 2000-2010 | Early stages of scholarship; beginning reception of Levinas |
| 2011-2015 | More explicit integration of ethics with social justice concerns |
| 2016-present | Development of multidimensional ethical frameworks |
| 2021-present | Stronger focus on combining philosophy with practical implementation |
| 2024-present | More empirical studies on ethics in mathematics education |

Figure 2: Timeline of research in ethics in mathematics education

---

[4] The connectivity graph was automatically generated using Litmaps on the full bibliography of this article. The tags are given by the last name of the first author.

Our analysis of the ethical turn in mathematics education will show that five distinguishing characteristics can be identified:

1. *Explicit Engagement with Ethics:* Authors identified with the ethical turn often draw upon philosophical ethics (e.g., Levinas) to theorize and frame their arguments, moving beyond observing and challenging sociopolitical realities to find new justifications for their underlying normative standpoints.

2. *Focus on Individual Interactions and Responsibility:* While group concerns remain important, the ethical turn provides an increased focus on the ethics of individual interactions, such as teacher-student relationships and the work of individual mathematicians.

3. *Questioning Mathematics Itself:* The ethical turn generally involves a renewed philosophical scrutiny of mathematics. Authors often challenge its perceived neutrality, purity, and universality (NPU) and examine the ethical responsibilities inherent in mathematical practice, including its applications and potential misuse.

4. *Bridging Discourses[5]:* Many authors surveyed here try to bridge gaps between different communities, such as mathematicians and mathematics educators. Framing concerns in ethical terms allows them to find new resonance spaces that, for example, critical theory could not. It can also involve mathematicians re-engaging with educational concerns out of ethical considerations for mathematics and its societal impact.

---

[5] The idea that ethics can serve to bridge competing discourses and different groups of actors (mathematicians, educators, etc.) was one of the motivations to write this analysis, to present it at such a length, and to not make use of a specific analytical perspective (e.g., a concrete version of discourse analysis). Scholars and mathematicians come to ethics in mathematics education with varying degrees of knowledge in the science and technology studies, history, and philosophy of mathematics. To make this review approachable was an explicit goal that shaped its structure, language and level of (philosophical) abstraction.

5. *Addressing Challenges of Poststructuralism:* Implicit in many authors is an engagement with challenges of the sociopolitical perspective, which Gutiérrez (2013, pp. 55–56) had already identified, such as concerns about relativism making "values and morals disappear" or prioritizing "intellectual exercise" over empowerment.

While perhaps not as substantive or widespread as the decades of sociopolitical scholarship in mathematics education, this emerging ethical focus adds a distinct dimension to the field. It prompts reflections on the ethical basis of teaching, the ethical responsibilities of educators and mathematicians, and the nature and impact of mathematics in the world, all by approaching it from a slightly different angle than before. We aim to explore the contours of this ethical turn and its relationship to the ongoing sociopolitical discourse within mathematics education.

## 1.1. Complications in understanding the ethical turn

The ethical turn, as suggested by Ernest (2024a) and analyzed in this article, needs to be contrasted with the ethical turn in the social sciences, where it has been characterized as a postmodern turn to the historically and culturally local that opens up a "shift to the discourse of the Others neglected in society and culture" (Canceran, 2023, p. 9). In mathematics education, this postmodern turn to the Other[6], and towards the marginalized, the excluded, and the value systems of diverse cultures, has already happened through the sociopolitical turn of mathematics education scholarship (see also Gutiérrez, 2013). For example, ethnomathematics presents an ever-present challenge of Eurocentrism (A. B. Powell & Frankenstein, 1997). The ethical turn takes on the theme of the Other, but, as we will show, the discourse changes its method and structure. The focus is now much more on the individual. The differences between the ethical turn in the social sciences and mathematics are subtle but noteworthy.

---

[6] In the discourse on ethics in mathematics education the Other (person) is often capitalized as a reference to Levinas; particularly the infinite respect that Levinas requests towards other people.

Furthermore, looking at the role of ethics in the social sciences, Lebow (2007) argues that ethics is fundamental to knowledge production, as it enables and maintains the set of practices necessary to deal with the pluralism of epistemological and methodological perspectives. From Lebow's perspective, it protects scholars from becoming polemics filled with intellectual intolerance. Such a perspective is also evident in the discourse on ethics in mathematics education, where some authors aim to use ethics to bridge differences, and others explicitly turn to the ethics of conducting mathematics education research.

However, there still remains an inherent difficulty and potential artificiality in trying to rigidly separate the "sociopolitical turn" from the "ethical turn" (e.g., as is evident in Buell and Piercey's (2024) discussion of ethical teaching in college classrooms). As such, we refrain from providing an exact definition that would be unnecessarily restrictive, but we still attempt to delineate the specifics of this new turn throughout our analysis. However, by focusing on an ethical *turn,* our reading of the literature may be slightly different from how others read this and the wider critical mathematics literature, and thus also what others have presented in their surveys, especially those who focused on the general positions that have or can be taken on (e.g., Dubbs, 2020).

A further complication is that debates on ethics in mathematics and its education are not limited to education journals but frequently occur outside the subject's core publication venues. As such, the following discussion also includes research that readers might not immediately identify as a piece on mathematics education. Nevertheless, as Gottesman (2016, pp. 167–170) already concluded in his historical review of critical pedagogy, there is tremendous value in reading broadly and outside the subject's core publications and engaging closely with the wider foundational literature underpinning more subject-specific discourse. The following analysis aims to follow this spirit, seeking to illuminate the discussions on ethics in mathematics education by attempting to find structure in an often explicitly poststructuralist discourse.

## 2. What constitutes ethics in mathematics education?

Ethics in mathematics education can be understood from multiple angles, including

1. the teaching of the ethics of mathematics,

2. ethics applied to mathematics education,

3. ethics applied to (general) education with a special focus on mathematics,

4. and mathematics applied to ethics in an educational context.

In the discourse on ethics in mathematics education, these perspectives can go hand in hand, both at the school and university levels. For example, authors rarely consider teaching the ethics of mathematics without considering the ethics of teaching itself. These various understandings of "ethics in mathematics education" mirror the general understandings of "philosophy of mathematics education," the latter of which Ernest (n.d.) explored in depth. Almost by definition, ethics in mathematics acts as a bridge between concerns regarding mathematics, education, and philosophy (Figure 3).

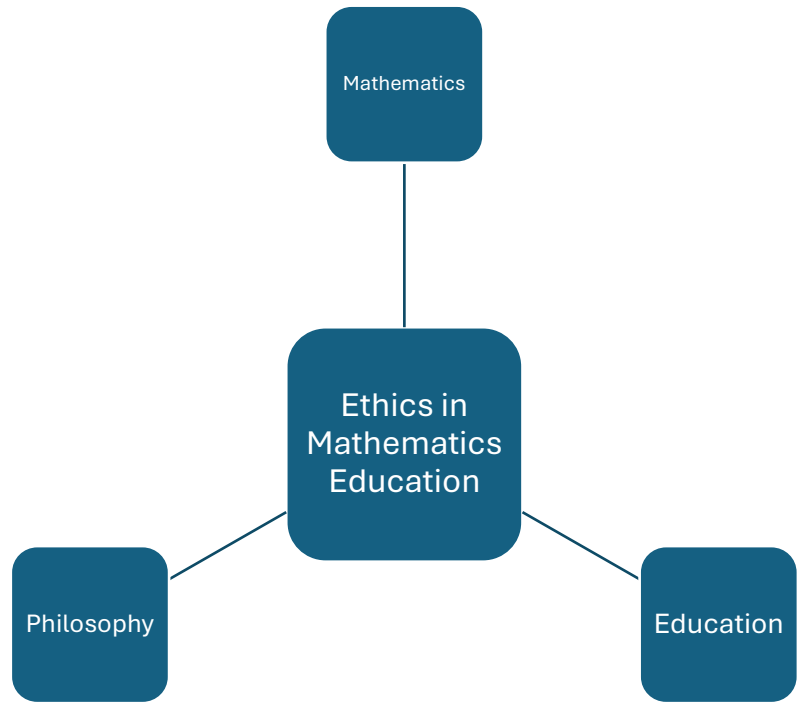

Figure 3: Ethics in mathematics education as a connecting force

## 2.1. Characterizations of Ethics in Mathematics Education

In the following, we will examine selective characterizations of ethics in mathematics education developed by participants in the discourse themselves. These characterizations will hint at the breadth of ethical questions asked.

Focusing largely on the context of school education, Boylan (2016) identifies four dimensions of ethics in mathematics education: the other person, societal and cultural problems of mathematics education, the ecological dimension of teaching and practicing mathematics, and the self as a practitioner, teacher or learner of mathematics. Taking the perspective of (research) mathematicians, Chiodo and Müller (2024a) identify concerns for the community of mathematicians (including learners of mathematics), concerns for the planet and wider socio-political world, and the development of mathematics as a field of knowledge. When compared, Boylan's approach to categorizing the dimensions of ethics in mathematics education is influenced by critical mathematics pedagogy. However, Chiodo and Müller write from a mathematician's perspective, merging the many "community concerns" that critical mathematical scholarship has tried hard to untangle, and they include an explicit ethical focus on mathematics itself.

Yet another approach is taken by Ernest (2019a, 2019b), who employs a common dichotomy of professional ethics. He separates the ethical responsibilities of the mathematics teacher into the duties of every human, every professional, and every mathematics teacher. Using a duty of care approach, he analyzes how five separate groups ("industrial trainers," "technological pragmatists," "old humanist mathematicians," "progressive educators," and "public educators") implicitly brought forth different ethical conceptions for the British mathematics curriculum in the 1980s and 1990s. These groups are mentioned here because we believe they still represent many stakeholders in mathematics education today. An orthogonal approach to structuring the ethical responsibilities of a mathematics teacher is taken up by Rycroft-Smith et al. (2024), who consider five levels of ethical, ranging from actively obstructing ethics in mathematics (level -1) to speaking out against bad mathematics (level 4). When we compare these two approaches, we see that the level-based approach is underspecified. What is considered ethical mathematics for the "Industrial trainers" in Ernest's analysis may not be considered ethical mathematics by progressive or public educators. While Rycroft-Smith et al. explore the depth of ethical awareness of individual actors, Ernest is looking at the breadth of possible positions that can be taken.

A further angle under which ethics in mathematics can be classified is opened up when one considers mathematical practices. Ernest (2021b) finds seven ethically relevant practices based on the distinction between "pure and applied, theoretical and practical, and with the mathematics utilized explicit or implicit" (ibid., pp. 7- 9):

1. "Pure theoretical practice with mathematics explicitly the main goal,"

2. "Applied theoretical practice with mathematics explicit,"

3. "Applied and theoretical practices with mathematics backgrounded,"

4. "Pure practical practices with mathematics explicit,"

5. "Pure and applied theoretical work, mathematics backgrounded,"

6. "Applied and practical practices, with mathematics explicit,"

7. "Pure and applied practical worth, with mathematics backgrounded."

These practices are translated by Bátkai (2023) into ethically relevant topics in the mathematics classroom: arithmetical errors, plausibility, faulty models and misinterpretations of models, text-based problems, questionable applications, automated decision-making, and modeling. When we compare Ernest's general classification with Bátkai's specialization in mathematics education, we notice that the foreground versus background roles that mathematics takes can be hard to untangle in the classroom and that the difference between pure and applied can reduce to a question of intent: Am I doing mathematics to enable me to do more mathematics or am I doing mathematics to solve a (real-world) problem?

These six classification systems of ethics in mathematics education show the complexity of understanding the ethical turn and how different authors can explore different concerns. One's perspective on ethics in mathematics education appears to be relative to one's positioning as an educator, the level of ethical awareness one has (or is willing to take), the area of ethical concern, how you're using mathematics and ethics, and most importantly the relationship between yourself and the Other. How these different aspects are related can differ among authors (compare Figure 4).

To illustrate the complexity, consider the following three hypothetical educators:

- Someone could explore level 0 (purity, neutrality and universality of mathematics) with a concern for the community and society, focusing on both applied and pure mathematics and coming at it from the perspective of a progressive public educator. Such a position would likely be found in ethnomathematics.

- Someone else might be concerned about the community of learners, conceptualized from an old humanist perspective, maybe even believing that there is (some) universality and neutrality in mathematics, with a strong focus on pure mathematics. Such a person would likely focus on the teaching of good mathematics.

- While another might be focused on taking a seat at the tables of power (level 3) with a concern for the community and a self-description of a progressive/public educator, reflecting strongly about the ethics of learning mathematics. This person may be focused on issues of diversity, equity, and identity.

These three variations show the immense complexity that can appear in scholarship on ethics in mathematics education, and this complexity has not yet explored the fundamental impact of the Self and the Other, the impact of culture, one's own individual arc in life, one's perspective on the nature of mathematics, or how all of these can change over the years. A fundamental relativism in the writings on ethics in mathematics education needs to be acknowledged.

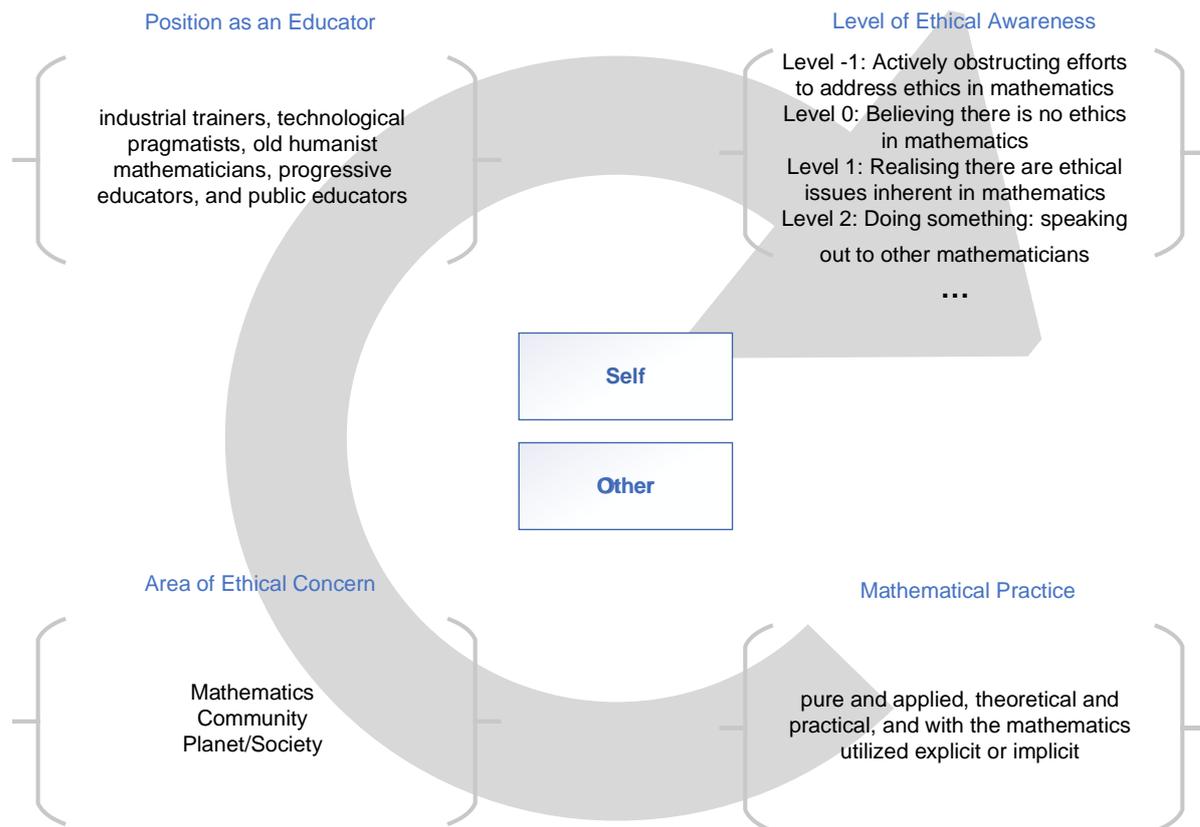

Figure 4: The complexity of ethics in mathematics education

This complexity also hints at a perspective that mathematics is a discourse itself rather than a fixed object to be practiced, taught, and learned, and so mathematics itself has immense formatting power, as Skovsmose (2023) also argues. The selection of authors, positions, and topics explored in this article tries to do justice to the depth and breadth of this complexity.

While it is certainly possible to structure an analysis along one or more of these classifications, this is not what we do. All these authors are actively involved in the discourse surrounding the ethical turn, and their classifications are necessarily biased towards their positions and objectives (e.g., Boylan writes from the perspective of education, Chiodo and Bursill-Hall write as mathematicians, and Ernest writes from the perspective of a philosopher). Taking on one of these classifications as the starting point would not grant us sufficient distance from the debate and likely mean we would miss something important. Instead, we decided to focus on narrative trends and argumentative lines, and thus, particularly on how various authors have recently started to use ethics in their work on mathematics education.

## 2.2. What is and is not part of the ethical turn?

Having explored some characterizations of ethics in mathematics education and explained what we aim to do, we can now explore what would and would not fall under the ethical turn. We showcase this through two examples.

First, and maybe least obvious, Bishop's (1994) research agenda on cultural conflicts in mathematics education does not fall under the ethical turn, at least not as conceptualized in this article. Building on ethnomathematical insights, his earlier works on mathematical enculturation (Bishop, 1988), and Keitel et al.'s (1989) UNESCO report on *mathematics, education and society,* Bishop structures his new research agenda around three narratives: the intended curriculum, its implementation in the classroom and the mathematics learned by the students. Looking closely at his proposed agenda, we can see that it is riddled with normative statements and pedagogical goals, but their ethical components are largely not made explicit. The author tries to unveil (and overcome) the sociopolitical issues and cultural conflicts in mathematics education. However, the central theme is conflicts and not ethics; although normative assumptions, statements, and goals play an implicit and explicit role in constructing this new research agenda, it is not the focus.

If we look at two works that would fall under the ethical turn, we can see a more explicit interaction with (philosophical) ethics. The following exposition will represent many authors in this analysis: they invoke ethics to support or question common assumptions and arguments in the literature, including critical mathematics pedagogy and ethnomathematics. Chen (2018) asks about the moral basis for teaching and learning mathematics, finding the justification for defining mathematics through pragmatist ethics and ethnomathematics rather than traditional arguments of the students' self-development, critical reasoning skills, or on other utilitarian grounds. Chen also challenges that social justice-themed mathematics can be the basis of such justifications, as it leaves largely "unquestioned the substance and the primacy of school math" (Chen, 2018, p. 155). While in agreement with Chen's overall criticism of contextless, socio-culturally unrooted mathematics, Mamlok (2018) nonetheless challenges her position, proposing that these conclusions cannot be deduced from and, in particular, using ethnomathematics alone.

The discourse between Chen and Mamlok highlights a common usage of ethics and one core challenge when trying to understand the existing mathematics education literature. Both authors attempt to find a new ethical basis for teaching and learning mathematics. In doing so, they build on but also have to leave the common sociopolitical arguments present in critical mathematics pedagogy. Furthermore, the authors do not subscribe to exactly the same understanding of critical mathematics, ethnomathematics, and their relationship. These nuances only come to light because the authors look outside sociopolitical discourses and extend existing debates on power, identity, and injustices by looking at what it means to act ethically. In the case of Chen, by first deriving an explicit pragmatist understanding of ethics, i.e., actions that are "communicable and conscious of the collective good" (ibid., 157), and only then trying to define the nature of mathematics using lessons from ethnomathematics. Chen attempts to disentangle the empirical, observable sociopolitical concerns, their ethical dimension, and the conception of mathematics precisely to find a way to bring these separate observations and thoughts back together in a more structured manner.

For many authors who are already part of the critical, justice-themed discourses in mathematics education, reflecting on (philosophical) ethics is a way to untangle the various dimensions of "critical" and "justice." Observable realities, normative standpoints, and the author's understanding of the literature can merge into a difficult-to-comb-out web of arguments that sometimes needs to be untangled. In short, for those authors, making ethics explicit is a way to make sense of existing problems and to give their (own) positions a new moral footing, which they could not find through sociopolitical considerations alone. In our reading, the authors are (consciously or subconsciously) trying to inject more structure into the variety of postmodern, poststructuralist arguments that make up much of critical mathematics pedagogy and other justice-oriented approaches to mathematics education.

In this spirit, Cannon (2018, p. 1331) asks what mathematics educators need to "think and unthink" to understand the role and potential of ethical thought in teaching and research, and Brown et al. (2016) go one step further. They ask about the role of theories and macro-trends in mathematics education in educational systems where "political exigency finesses educational ethics" (ibid., 288). Both ask us to reflect on the role that critical perspectives (should) play in mathematics education. The ethical turn in mathematics and its education is precisely such a reflection; it leads to a nuanced extension of existing sociopolitical scholarship and views of the nature of mathematics.

## 2.3. Doing good versus preventing harm

There are two partially overlapping themes in the discourse: ethics in mathematics education as a force to do good or as a way to prevent harm. For many authors, the questions of doing good and preventing harm are deeply connected, but some explicitly focus on doing good by preventing harm from mathematics, while others primarily seek to do good with their mathematics. For example, based on her experience of training Indigenous mathematics teachers, Abtahi (2022b) asks, "What if [she] was harmful [by teaching Western mathematics]"? Similarly, Ernest (2018) poses the question of whether mathematics is harmful to its learners and examines the "collateral damage" caused by learning mathematics (Ernest, 2016a), leading him to consider aspects of medical ethics with a distinction between no-maleficence and beneficence (Ernest, 2020). As analyzed by Müller (2024) and discussed later in this article, the literature from the (former Cambridge University) Ethics in Mathematics Project is also largely focused on preventing societal harm. Skovsmose (2021) also considers the potential harm from mathematics by looking at the potential relationships between mathematics and a crisis: one can use mathematics to get a "picture [of] a crisis", mathematics can be a "constitutive [part of] a crisis", and mathematics can "format a crisis" and thus the potential actions taken (ibid., 369).

Many other scholars, including those fostering democratic participation, regularly focus on doing good and empowering students (cf. Aguilar & Zavaleta, 2012; Buell & Shulman, 2019; Skovsmose & Valero, 2013). A concrete example is given by "mathem-ethics", whereby geometry is used to teach ethics, emotions, and self-development by comparing geometric figures and their loci to similar emotional states at school and in prison (Cerruto & Ferrarello, 2023; Ferrarello & Mammana, 2022).

The duality of preventing harm and doing good is also taken up by the various calls for and considerations of Hippocratic oaths in mathematics and its education. Building on the concept of ethical filtration (Skovsmose, 2008, p. 167), i.e., how ethical nuances are reduced and lost when real-life problems are framed in mathematical language, Freitas (2008) worries that curricula focusing on real-life problems still prioritize mathematical skills and knowledge over social justice concerns. According to Freitas, an oath would be beneficial to promote ethical reflection in an otherwise largely morally unreflective field, and such induced ethical reflection in the classroom would necessarily lead to improvements in social justice as a "source of our sense of agency and action" (ibid., p. 88). Seeing the power that mathematics wields in modern societies, Freitas (2008, p. 92) further asks, "Why wouldn't [...] an [Hippocratic] oath accompany the application of any powerful tool to the world?"

Rittberg (2023a) and Müller et al. (2022) analyze various other calls for Hippocratic oaths by (research) mathematicians and educators, finding their potential necessity to instigate and guide ethical discourse, but also their insufficiency to solve the ethical questions of mathematics and its education. The institutional structure, value systems, and professional identity among mathematicians are seen as insufficient compared to medicine. Their analyses are suggestive of the tendency that existing calls for Hippocratic oaths have been largely ignored within the wider community, likely either because of the neutrality, purity, and universality of the mainstream[7] view or because of a postmodern/-structuralist rejection of universalism. Here, an oath represents a

---

[7] We decided not to define the mainstream view of mathematics because we think that one generally has a fairly good idea of what it means, and because even small surveys among mathematicians show that it is hard to pin down exactly (e.g., Sriraman (2004)), as

 "[t]he working mathematician is a Platonist on weekdays, a formalist on weekends. On weekdays, when doing mathematics, he's a Platonist, convinced he's dealing with an objective reality whose properties he's trying to determine. On weekends, if challenged to give a philosophical account of the reality, it's easiest to pretend he doesn't believe it. He plays formalist, and pretends mathematics is a meaningless game." (Hersh (1997, p. 197), as quoted in Greer et al. (2024, p. 2)) universalist approach to ethics that is then interpreted as incompatible with poststructuralist perspectives.

### 3. A common "starting point"? Neutrality, purity, and universality (NPU) of mathematics

A central concern within the ethical turn in mathematics revolves around the contested notions of mathematical neutrality, purity, and universality (NPU), which "coagulated as a paradigmatic dogma portraying science and mathematics as isolated phenomena not involved in any social or political complexities" (Skovsmose, 2024a, p. 13). We use the term "purity" rather than "objectivity" to signal the ethical dimension of seeing mathematics as objective. Purity has a double meaning of *not being mixed with something* (in the case of mathematics, understood as consisting of only rational objectivity) and *being free from harm* (and thus, by definition, ethically positive). In our experience and after reading the literature, debates about the ethics of objectivity are often so deeply connected to the stipulated purity of mathematics that, in our eyes, it makes sense to focus on it immediately.

Similarly, we merge neutrality, purity, and universality. Even though some authors only speak about one, the other two always appear to be lurking in the background; universality or purity can come as an argument for neutrality (If something is universal or purely rational, how can it not be neutral?) etc. Each of the three aspects already has the capacity to make an open discourse on the ethical nature of mathematics impossible: the case of neutrality is explored in-depth by Ravn and Skovsmose (2019), the case of purity is explored by Müller (2018a), and the case against universality is the cornerstone of ethnomathematics.

## 3.1. Why do authors challenge NPU?

For different authors, this concern comes at different points in their chain of reasoning. For scholars already deeply engaged with critical mathematics education, particularly those building on the vast trove of postmodern scholarship and beyond, challenging the NPU of mathematics can be (deductively) given by the rules of said scholarship. Universality must be challenged; plurality is the norm in their discourses, and the "incredulity"[8] towards the meta-narrative of PNU follows (almost immediately) for many. However, for mathematicians joining the ethics in mathematics discourse, ethical engagement with PNU often takes more of an abductive line of reasoning. They may have encountered an ethical issue within (their) mathematics and now attempt to trace it back to its origins. And yet others may have long seen the issues and are now trying to construct new rules of the game, i.e., rules for ethical reasoning about mathematics and practicing mathematics, or argue mathematics itself is an essential driving force behind today's plurality of truths (e.g., Rottoli, 1998). The latter positions may posit that mathematicians should be less concerned with finding ultimate truths but reflect more on their freedom within mathematics and, thus, their ethical responsibility to deal with that freedom.

In short, how individual authors come to challenge the NPU of mathematics varies across disciplines, as well as individual experiences with mathematics and their conceptualization of the nature of mathematics. A careful reading of authors is necessary to understand why and how they are challenging the neutrality, purity, or universality of mathematics. Such care is necessary when interacting with the literature and in person. As is explored in Müller (2024, p. 72), these differences can hinder productive discussion if not articulated properly. We see that the complexity of ethics in mathematics comes not just from the variety of potential positions and characterizations we explored earlier but that the origins of one's thinking can be of similar importance. These may or may not be correlated.

---

[8] In reference to Lyotard (1984).

Furthermore, when education scholars and others begin talking about the socio-constructive nature of mathematics, it can be read by professional mathematicians as an attack on their discipline and their hard-earned mathematical truths (for a general discussion of this phenomenon in the sciences, see Vessuri (2002, p. 136)). For example, Ernest (2024c, p. 1241) argues that "attempts to raise ethical issues with regard to pure mathematics are seen as possibly tainting or lowering the subject from its elevated state of purity." In a personal anecdote, Maaß (2020) similarly states that his attempts to introduce ethics into discussions about mathematics and its education were met with "disinterest" (at best) and "warnings" (at worst), providing the anecdotal quote to "leave the topic alone, do not get yourself in trouble" (ibid., 2020, pp. 50–51; authors' translation from German). Gordon (2022, 2025) considers a communication ethics-specific lens through which they analyze the ethical dimension of the language of mathematics; in particular, finding that its focus on presenting finished results and its inherent emphasis on excluding investigative language can lead to the exclusion of students and hinder their creative processes.

The tension between the "commitment to an aesthetic of concision promulgated by tradition" and excluding the "investigation that informed the many steps of decision-making" (ibid.., 29) lies at the heart of his ethical analysis. This tension between mathematics and its practice – sometimes implicit, sometimes explicit – is often present. These streams of conflict can also be neatly explained by the statement that "[mathematicians] have claimed to be useful mathematicians doing useless mathematics [when they are confronted with questions of an ethical nature]" (Karaali, 2019, p. 12), a self-description in stark contrast to the increasing "power differential" (Buell & Piercey, 2022, p. 3) between mathematicians and non-mathematicians. In short, there are different reasons and motivations for challenging and defending neutrality, purity, and universality.

However, questions regarding the NPU of mathematics are crucial for those participating in the ethical turn. For some authors, it even represents the lowest development of personal ethical consciousness and a stumbling block to overcome (e.g., Chiodo & Bursill-Hall, 2018; Chiodo & Müller, 2018), while others think about the issue from a histo-philosophical perspective and have begun to look for philosophies that allows them to conclude an ethical responsibility in the presence of mathematical purism (e.g., Ernest's (2021c) discussion of MacIntyre's virtue ethics). These two approaches represent two extremes: how much of the perceived NPU of mathematics constitutes a deficit to be overcome, and how much of it can we tolerate while still being able to discuss ethics in mathematics and its education?

In general, discourses coming out of the corner of critical mathematics pedagogy are immensely critical of any claims to universality and, thus, the neutrality and purity of mathematics (Müller, 2024). However, the level-based approach developed by Chiodo and Bursill-Hall shows that pragmatic approaches can also follow a similar trajectory, even when they attempt not to subscribe to a particular educational or philosophical line of thought. The latter happens primarily because such pragmatic approaches are still fundamentally postmodern: the pragmatism from "not wanting to take a position" acknowledges the plurality advocated for by critical mathematics pedagogy.

## 3.2.   How do authors challenge NPU?

For mathematics education, the NPU debate has immediate implications. For many authors, it influences — and sometimes determines — how much ethics can (and should) be incorporated into the classroom. A teacher subscribing to certain neutrality assumptions may be less likely to consider other moral issues in the classroom, mathematical institutions may push back against social justice initiatives (e.g., Crowell, 2022), and governing bodies are more likely to introduce curricula that leave little space for ethics (as discussed in Ernest, 2019a). Of course, this theme has already been present in sociopolitical works on mathematics education, and it continues to play a central role in debates on ethics in mathematics. How much ethics there can be in the classroom is fundamentally also a question of how much universalism and neutrality of mathematics exist in the minds of those involved. In short, sociopolitical discussions often lead to deeply connected ethical and metaphysical questions. The most direct case of this is probably when Ernest (2012) asks if ethics can serve as the first philosophy of mathematics education rather than epistemology or others.

The question of the NPU of mathematics also arises when authors discuss the ethics of mathematical applications. For example, Straehler-Pohl (2016) builds on Morozov's (2013) analysis that the increased quantification and algorithmizing of life does not happen despite ethics but precisely because mathematicians, engineers, managers, and others involved in it are working for what they perceive to be a better world; a world where data, mathematics, and technology can and should be the solution to any, and all problems. Morozov calls this ideology solutionism, and Straehler-Pohl uses it to contrast it with critical mathematical pedagogy. By juxtaposing the radically different perspectives on good human life and a better world, this opposition quickly evolves into a discussion on the theoretical underpinnings, where critical theory (underpinning critical mathematics pedagogy) is opposed to the NPU of mathematics (underpinning solutionism). The analysis can be read as a call to rethink how we present school mathematics and the role that mathematical solutions play in it. Here, we see how questions about the NPU of mathematics can quickly bridge the gap between the classroom and the outside world.

Radford (2008) challenges philosophies of mathematics that understand the subject as an objective, individualistic, and purely rational endeavor. He suggests that all learning is cultural learning, achieved through interaction with cultural artifacts and social communities, which extends beyond acquiring knowledge to "becoming someone." Radford particularly stresses the social, cultural, and historical dimensions of mathematics and learning, whereby learners of mathematics do not just get to know mathematics but find their place within the socio-histo-cultural context of mathematics. Ethics in mathematics education is thus central, and it largely corresponds to the teacher's responsibility in shaping students' mathematical identity and their position in and relationship to culture.

### 3.3. What do philosophers say?

Philosophers of mathematics also discuss the moral challenges that arise from and are inherent in the practice and teaching of mathematics. For example, Wagner (2023) relates ethical questions to the development of today's formal approach to mathematics, developing an approach informed by social philosophies of mathematics rather than critical mathematics pedagogy. In turn, he finds a different angle of critique by explicitly pinpointing Hilbert's formalistic approach as ethically problematic:

"Hilbert's victory (not his dream of a formal proof of consistency, but the successful banishment of philosophical debates from mathematics in favor of a pluralism of axioms and formal systems) guaranteed a renewed consensus regarding the validity of mathematical proofs relative to well-specified systems, breathing new life into the image of mathematics as a domain of universal certainty. This was an ethically problematic achievement, as it depended, to a large extent, on exiling from the realm of mathematics all disagreements over the validity of arguments that could not be settled by Hilbert's formalizing […] [which] in turn, had a lasting impact on the over-valuation of mathematics […]." (Wagner, 2023, p. 14 building on Wagner, 2022).

Other philosophers have explored the ethical dimension of participating in mathematics through the lens of epistemic injustice, considering when inclusion (e.g., of learners) and exclusion (e.g., of "mathematical cranks") is warranted, paying particular attention to aspects of common knowledge (Rittberg et al., 2020), the apprenticeship model of learning mathematics (Tanswell & Rittberg, 2020), and the variety of reasons that can lead to exclusion on epistemic grounds (Rittberg, 2023b).

Another assessment of the ethical questions surrounding pure mathematics is presented when Hersh discusses that "the ethical component is very small. Not zero, but so small it's hard to take very seriously. […] That's why people say pure mathematicians live in an ivory tower. One answer to this could be, 'Well this is fine! There is no need for mathematicians to have a code of ethics because what we do matters so little that we can do whatever we like.' […] I find it rather scary. Because it means that if we become totally immersed in research on pure mathematics, we can enter a mental state which is rather inhuman, rather totally cut off from humanity. That's a thing we could worry about a little bit." (Hersh, 1990, p. 22) His way of thinking raises another concern: rather than asking if mathematics *is* disconnected from humanity, he asks if good mathematics *should* be.

Questions about the NPU of mathematics are also connected to those surrounding the differences and similarities between moral and mathematical reasoning and how they can be at odds with each other (e.g., Clarke-Doane, 2015). Such a question is particularly employed by Nickel, building on Kambartel (1972). The latter traces the historically dominant separation of ethics and pure mathematics to Max Weber's call for a science free from values. Discussing his contribution to the field, Nickel (2022) shows that whether reasoning in ethics and mathematics shows parallels (and are indeed not separated) has historically depended on whether ethics and mathematics were understood as rational or irrational, finding positions for and against separation in either case (also discussed in Nickel, 2015). Nickel (2005) finds important similarities, including statements of universality, plenty of historical contact between the two disciplines (including the search for (eternal) forms that mathematics and ethics had in common for a long time), and structural parallels (e.g., freedom and rules). Two distinctive differences are the role of individual exceptions and the other (person/mathematician) (ibid., 428). For Nickel, studying similarities and differences in mathematical and ethical reasoning can bridge the gap between ethics and mathematics.

### 3.4.   NPU, institutions, and teaching

We see here how some begin their arguments using empirical observations (e.g., observable injustices) or arguments from the critical social sciences (e.g., the power that mathematics holds in society), while others immediately return to the philosophical foundations of ethics to challenge the NPU of mathematics. Either line of argumentation has renewed interest in questions about the philosophical nature of mathematics and the role of mathematics in society, including the role of mathematical institutions. Pearson, then Executive Director of the Mathematical Association of America, observed the contrast between the codes of ethics of the American Mathematical Society and the emerging literature on ethics in mathematics. He argued that major mathematical societies too often still understand mathematics as ethically neutral. He explored the tension between the majority view of NPU and the corresponding debate within the ethics in mathematics literature through three blog posts (Pearson, 2019a, 2019b, 2019c). Pearson concludes the series with an argument about why a "concerned, quantitatively literate citizenry" requires us to learn more about questions like these.

As this discussion shows, questions about NPU are often a mix of ethical, sociopolitical, and wider philosophical concerns, and for many education scholars, they come together with calls for reform. Such calls for reform, however, are themselves noteworthy of ethical analysis, as Grootenboer (2006) puts it:

"Given the many calls for belief reform with preservice primary teachers, it is also likely that mathematics educators will be actively and overtly promoting the views that they consider 'correct' and appropriate. I am not suggesting that it is necessarily wrong for mathematics educators to do this, but it is an ethical issue when someone is trying to change the beliefs of another, particularly when the educator also has a powerful gate-keeping role. I think it also behoves mathematics educators to carefully consider the beliefs they are promoting and to thoroughly research and analyse the benefits of these views for effective mathematics education." (ibid., p. 274)

In summary, debates about the neutrality, purity, and universality of mathematics should not always be understood as the starting point of one's reasoning ethics in mathematics. However, they constitute an important intermediate step in the reasoning of the involved scholars, often after the scholars have done their "archeology"[9] of its origins. Furthermore, these debates do not always happen on neutral grounds, but they often come with specific educational goals in mind. To at least find some common ground in the "NPU debate," educators and mathematicians might find it useful to use the related idea that mathematics is indefinite with respect to concepts, proofs, topics, applications, cultures, and power (Skovsmose, 2024b), as it is not yet as value-laden.

## 4. How much, and what philosophy as the foundation?

When the reviewed scholars engage with ethical questions surrounding mathematics education, they use various philosophical traditions to ground their arguments and justify their normative position. The choice of philosophical foundation varies, reflecting the diverse ways ethical issues are conceptualized within the field. Some researchers draw on established ethical theories, while others look at specific thinkers to illuminate the moral dimensions of mathematics teaching, learning, and research. Nonetheless, common themes can still be identified.

Radford (2023) argues that ethics is inherently existent in mathematics education, as any form of teaching and learning always involves interactions between the teacher and the learner, and thus sociopolitical issues, including the validation of specific forms of teaching, learning, and mathematical knowledge. He analyses the different influences that Hobbes and Kant have had on our understanding of mathematics and the classroom: Hobbes' contractualism has influenced transmissive learning philosophies, while Kant influenced constructivism. In the spirit of postmodernism, Radford leaves such scholars behind and argues for an ethics based on mutual responsibility, care, and commitment to each other to foster social justice-themed issues such as inclusion and respect. We will see that such a transition is fairly common among scholars. Their selection involves philosophical positions that challenge the mainstream NPU view of mathematics and focus on the responsibility of the individuals. In the following, we will discuss three more common approaches before reflecting on minor trends.

---

[9] Not rarely does this digging happen with a traceable line to Foucault (1970), even leading to explicitly Foucauldian perspectives on ethics (e.g., Dubbs (2024)).

The three important trends are the "ethics as a first philosophy"-approach by Emmanuel Levinas, the selection of non-Western scholars with departing viewpoints, and approaches building on (philosophical) pragmatism. We picked these not just because of their relevance to the discourse but also because of their diversity: the strand concerning Levinas is particularly focused on a single philosopher, the non-Western philosophies represent the transnational and intercultural aspect of ethics in mathematics education, and pragmatism is rather more applied and sometimes seen as less philosophically deep. Together, these three provide us with sufficient breadth and depth to get a glimpse into how individual characteristics of the discourse have developed over time.

## 4.1. Reception of Emmanuel Levinas

Building on a critique of the mainstream view of mathematics as universal, neutral, and pure, scholars frequently seek alternative groundings to Emmanuel Levinas's work. Having problematized the (modern) perspectives of mathematics that often obscure or negate differences, these researchers find a new point of departure in Levinas. His philosophy prioritizes the ethical encounter with the irreducible Other before questions of knowledge or being. This perspective makes Emmanuel Levinas (1987) a widely received philosopher in the literature on ethics in mathematics education. His approach to rejecting traditional Western philosophies in order to prioritize ethics as the first philosophy, combined with his analysis of the inescapable 'Other' (person) and the corresponding obligation to treat the Other with the utmost respect and moral behavior, is well-received among scholars focused on interactions in a mathematics classroom. Moreover, Levinas's often radical language and promotion of the Other's inescapable self-worth align with the normative claims made by critical mathematics pedagogy and similar approaches to social justice-centered mathematics education, particularly his request "to have to answer for one's right to be, not with reference to the abstraction of some anonymous law, some legal entity, but out of fear for the other"(Levinas, 2022, p. 27) and his questioning if "my 'place in the sun' [...] [was] not already a primal usurpation of places that belong to the other, already oppressed or hungry before me"(ibid., p. 27). He appears as a suitable philosopher for critical scholars who, focusing on a poststructuralist and postmodern rejection of universal moral claims of rational thought, aim to rethink mathematics and its education.

The first reception of Levinas on issues of ethics in mathematics education was likely given by Neyland (2002) and his criticism of neoliberal management thought in (mathematics) education. In his view, the latter damages the ethical relationship between students and teachers. Neyland uses Levinas to move classroom interactions away from technical mathematical knowledge to focus on the ethical interactions between individuals in the classroom. Building on this insight, Atweh's work on "socially response-able" mathematics education represents the most extensive usage of Levinas in the discourse on ethics in mathematics education. He criticizes mathematics education as overly focused on developing the rational mind (Atweh, 2007) and argues for the training of mathematics teachers that enables them to respond to ethical questions in the classroom ("response-able"), in particular, the ability to react to the (welfare) needs of the students (Atweh & Brady, 2009). A case study of three teachers displays how "beliefs about the nature of mathematics, its role in the curriculum and the real world determine [a teacher's] readiness to take risks in changing classroom practices"(Atweh & Ala´i, 2012, p. 104), highlighting the significant effect that a perceived purity and neutrality of mathematics has on ethical classroom interactions.

Building these findings, Atweh (2013b) further extends his engagement with Levinas, criticizing attempts to codify ethics (in mathematics education) as they, for him, generally follow a problematic (Western) dominance of reason and rationality. His solution is again to reintroduce Levinas's moral philosophy into the discourse, as his radical acceptance of the Other works well for the sociopolitical scholarship surrounding him. This engagement with Levinas is particularly social justice-centered, focusing on empowering students' self-development in all aspects of life and social participation (see also Atweh, 2012, 2013a). Issues in the construction of equity and identity further lead Atweh and Swanson (2016) to discuss the diversity of intellectual approaches to equity, using Levinas' focus on being in a relationship with the Other to shift the discussion to the normative circumstances underpinning the discourse. In doing so, they attempt to gain new ethical insight into those learners who are left behind, using Levinas's ideas to justify a moral obligation to equity.

Similarly, Maheux and Roth (2012) argue for an approach to mathematics education grounded in the ethics of Levinas, discussing co-teaching and cogenerative dialogue as tools to encourage responsible, mutually shared learning experiences focused on the students' and teachers' benefits. Maheux and Proulx's (2017) engagement with Levinas leads them to focus on collaborative interactions in the mathematics classroom, focused on humanistic values rather than deficit-based banking approaches to teaching. Again, these works focus on accepting the Other in all classroom interactions.

The theme of ethical classroom interactions is also taken up by Demattè (2022b), who analyses how mathematics teaching creates situations of ethical importance, particularly when it hinders the students' ability for self-development and self-expression, which is often represented by an inability to make sense of the mathematical content they are being presented. He understands these teacher-student interactions as subtle forms of violence and extends the arguments to student-text interactions, whereby the interaction of a reader with a mathematical text is also (indirectly) a meeting between the reader and author (also explored in Demattè, 2022a). In this analysis, Dermattè pays close attention to the exclusion of students by a group or community of insiders with mathematical knowledge (including teachers and professional mathematicians). He argues that ethical mathematical texts must balance mathematical content and educational choices in ways that respect students as their readers. Like scholars before him, Demattè finds the moral framework to resolve these issues in Levinas and his radical acceptance of the Other.

An alternative but well-aligned perspective is also given in Paul Ernest's (2012) quest for a first philosophy of mathematics education, where he follows Levinas' conclusion. Just like Levinas, Ernest finds such a first philosophy in ethics, noting in particular, that grounding the philosophy of mathematics education on ethics has "immediate implications for individualistic learning theories" (ibid., p. 14). For Ernest, the encounters between different people, including teachers and students, in a learning setting must henceforth be ethically understood before specific teaching practices and educational philosophies are constructed.

The overall reception of Levinas in the wider educational sciences has followed different goals and discursive trajectories (Zhao, 2019), but he is overall well-received by (critical) scholars as he "has addressed many problems of modern Western thinking" (Zhao, 2016, p. 1) that also trouble (critical) educational scholars, in particular, challenging the "quest for certainty, uniformity, and accountability" (ibid., p. 4) which has been brought to many education systems. However, despite a more varied reception in education, his reception in mathematics education focuses mainly on classroom interactions (teacher-student, student-student, reader-author) and the rejection of totalitarian tendencies of rational thought. Levinas's reception is strongly influenced by existing poststructuralist discourses of alterity, identity, and power in critical mathematics education, including what Walshaw (2013, p. 116) calls the "widespread concern for 'outsiders' excluded from the mainstream story of mathematics." These pre-existing discussions laid the foundation for making his radical acceptance of the Other easily digestible for critical mathematics educators looking for new justifications for their normative standpoints.

Unlike the complex reception of Lacan's psychoanalytic study of the Other in critical pedagogy (e.g., Armonda, 2022) and mathematics education (e.g., Moore, 2024), the usage of Levinas by (critical) mathematics educators is so far fairly uniform in its usage and goals, i.e., the normative justification of caring for the Other and a justification for putting ethics first in the philosophy of mathematics education. This uniform reception may be in part because even though Lacan's Other is still there before us (Belsey, 2002, p. 58), he still "seeks a subject who frees themselves from the bondage and dominance of the desire of the [O]ther" (Karimi & Asghari, 2021, p. 133), requiring a more intricate usage of his ideas in a discourse that is already very much centered around the concern for the Other. The constant struggle in the ethics in mathematics literature with the (perceived) neutrality, purity, and universality of mathematics, and thus the constant rubbing against traditional Western philosophies, has been empirically explored by Atweh and Ala´i (2012) when they study their effect of teachers' preconception of mathematics. Their findings and the overall (theoretical) discourse on Levinas suggest that (early) training future mathematics teachers may be beneficial to counter the NPU problem and allow teachers to focus more on the ethics inside their classroom.

---

The lesson here is that the acceptance of the Other[10], i.e., the other learners, teachers, and mathematicians, needs to triumph over potentially harmful mainstream ways of thinking about and teaching mathematics. Humans are at the center of mathematics, but there should not be anthropocentric mathematics.[11]

---

## 4.2. Non-Western philosophies

The dependency of mathematics on history, culture, and sociopolitical contexts is a cornerstone of the sociopolitical turn in mathematics education and in the analysis of mathematics as human practices, with far-reaching consequences for social justice approaches (Mukhopadhyay & Roth, 2012). Furthermore, looking at mathematics from a non-Western perspective has a long tradition, particularly in, but not restricted to, ethnomathematics.

Out of ethnomathematics – especially the understanding that there is not one, but many mathematical cultures and ways of thinking mathematically (Albanese et al., 2017) – new approaches to teaching mathematics have developed, often focused on finding a balance between teaching dominant mainstream views of mathematics and non-dominant perspectives (e.g., Peralta, 2020). Such critiques of Eurocentric discourses on mathematics education, building on tensions between Indigenous forms of knowledge and the curriculum (e.g., Sterenberg, 2013), hopes to decolonize mathematics (e.g., Lunney Borden, 2021), and the lived identities and cultures of students and teachers (e.g., Luitel & Taylor, 2007), have brought forth a plethora of what can be called non-Western critiques of mathematics and its education. Such works generally give a vast perspective of the privilege that the (Western) mainstream view of mathematics has in curricula around the world and its "culturally homogenizing force"(Namukasa, 2004, p. 209).

---

[10] The "Self and the Other" in the mathematics classroom has also been explored by Radford (2024), connecting it with ethics of authenticity as developed by Taylor (2003).

[11] The question of how to develop an ethics in mathematics that fosters multispecies flourishing has also been explored by Khan (2020), writing: "The short thesis of this paper is that current conceptualisations of a 'Mathematics for Human Flourishing' (2020), while a welcome and important contribution to ongoing conversations in our fields, continue in a tradition of an anthropocentric mathematics that fails to attend to the evidence that there has never existed widespread human flourishing anywhere without simultaneous attention to multispecies' flourishing" (ibid., p. 231). The idea of a mathematics entirely focused on living without killing, i.e., for the survival of all and everyone, has also been explored by D'Ambrosio (2024).

These approaches are inherently focused on ethics, often explicit, sometimes implicit, as D'Ambrosio (2006, p. 1) writes: "Ethnomathematics is embedded in ethics, focused on the recovery of the cultural dignity of human being." However, unlike Ernest's (2024a) view of the ethical turn as strongly focusing on individual interactions, the concerns of ethnomathematics are generally focused on cultural and group ethics. It is an "ethics of diversity [focused on the] respect for the other (the different); solidarity with the other; cooperation with the other" (D'Ambrosio & D'Ambrosio, 2013, p. 22). For example, D'Ambrosio and Rosa (2017, p. 287) give two goals for bringing ethnomathematics into school curricula:

1. to "demystify school mathematics as a final, permanent, absolute, unique form of knowledge"
2. and to "Illustrate intellectual achievement of various civilizations, cultures, peoples, professions, and genders."

However, some scholars also consider the ethics of individual learners, teachers, and mathematicians. One such example can be found in the works of Abtahi (2021, 2022a, 2022b). By connecting the quote, "The universe is not outside of you. Look inside yourself; everything you want, you already are." by the Persian philosopher Rumi and the Anishinaabe Indigenous saying, "If you seek your vision around the circle in deeply thoughtful and respectful ways, you will find something that was previously hidden to you," Abtahi (2022a) bridges the world between her (Canadian) Indigenous students and her Persian heritage, to find an ethical perspective of her actions in the classroom. She asks what it means for a "learner and teacher of mathematics […] to take a journey into [themselves]? And why does it matter?" (ibid., p. 1), arguing that understanding the self is central for the ethical mathematics teacher.

Nevertheless, as this example shows, sociopolitical and ethical issues are near-impossible to untangle in such a context, especially when scholars begin their journey into ethics in mathematics education from personal experiences in the classroom. In such cases, identity, politics, and culture are necessarily interwoven with ethical considerations. Nonetheless, we can see the bridging capabilities that ethics can have; it can establish a bridge between a (culturally sensitive) educator and their students, between dominant and non-dominant views of mathematics, and between mathematicians and potential users or learners of their mathematics. As the close connection between D'Ambrosio's ethics of diversity and Abtahi's experience highlights, the power of considering non-Western philosophies does not just lie in the fact that they are, in a certain sense, different from the mainstream view of mathematics but that they can fulfill important identity functions without which mathematics learning as cultural learning (as explored in Radford, 2008) cannot happen.

The lesson here is that mathematics teachers need to be able to deal with the ambiguity surrounding mathematics and their teaching. Teachers need to be able to navigate Nepantla, the "uncomfortable space where there is no solid ground, that has no official recognition" (Gutierrez, 2012, p. 35), and to develop an ethical and political identity to do so (Gutiérrez et al., 2024). They need to self-critically understand the ambiguous space left between themselves, their morality, and the politics surrounding them to accept others and challenge unethical mathematics in the many interactions and relationships that develop in the classroom.

## 4.3. Pragmatic approaches

The pragmatic school of ethics in mathematics appears to be primarily advocated for by the (Cambridge University) Ethics in Mathematics Project (2018). The project mostly focuses on the ethics of mathematics at the university level. It had an initial (yet recurring) focus on the question of the neutrality of (pure) mathematics (Chiodo, 2020; Chiodo & Clifton, 2019; Chiodo & Müller, 2024b; Müller, 2018a) and the various "levels of ethical engagement"[12] (Chiodo & Bursill-Hall, 2018; Müller, 2018b) that mathematicians can display in their work. These levels were later expanded to an ethics framework for mathematics teachers, explaining the various roles and positions they can take in an ethical discourse at school, now including "actively obstructing efforts to address ethics in mathematics" while still ranging to "calling out the bad mathematics of others" (Rycroft-Smith et al., 2024, p. 377).

---

[12] How to bring a mathematician from believing that mathematics is neutral to an understanding that there is ethics in mathematics has been explored in many ways; one interesting approach is given in the fictional dialogues written by Ernest (2016b, 2021a).

For the latter, Koblitz's (1981) criticism of the misuse of mathematics in the social sciences or the publication of a statement by the Royal Statistical Society (Inference, n.d.) condemning the misuse of statistics in court would be two examples where mathematicians called out bad mathematics, and Shulman's (2002) discussion of ethically questionable exercises in mathematics problems would be an example related to the mathematics classroom. Leinster's (2014) earlier discussion of how mathematicians need to realize that their work can be misused is a good representation of difficulties surrounding level 0, i.e., realizing that there is ethics in mathematical practice.

Inspired by the role that mathematical modelers played in the United Kingdom's fight against the Covid-19 pandemic, the project developed an approach of 9 (Chiodo & Müller, 2020) and later 10 pillars of ethical mathematical practice (Chiodo & Müller, 2023). The pragmatism of the approach is displayed in this manifesto, which, following a prototypical mathematical workflow, lays out a sequence of (in-depth) moral questions for the mathematician to answer. Noteworthy is that, unlike other approaches discussed in this survey, the approach is solely focused on applied mathematicians, does not attempt to provide answers but questions for mathematicians to ask, and explicitly tries to avoid taking on a specific normative standpoint. However, a closer inspection reveals various parts of philosophy that have influenced the project, including AI ethics, diversity studies, computer science, engineering ethics, and debates on trust and the political nature of technological artifacts. The manifesto was later expanded to set the foundation for an approach to embedded ethics in mathematics, i.e., the explicit teaching of ethics in mathematics in mathematics courses (Chiodo et al., 2023) and arguments for the political nature of mathematical artifacts (Müller & Chiodo, 2023).

In large part, the project does not ascribe to all the underpinnings of critical mathematics pedagogy(Müller, 2024). However, it tries to focus on bringing applied ethics to a potentially unwilling group of mathematicians and mathematics students who ascribe to philosophies of mathematics that posit its purity, neutrality, and universality. Additionally, the project appears to be less concerned with student well-being and mathematical achievement than with the potential harm mathematics can do in society. The focus is outward, and much of the literature is written from a position of privilege[13] that does not have to worry about students' mathematical achievement.

Other authors have also explored ethics in mathematics education more pragmatically. For example, Porter (2022) offers a practical reading of data ethics for mathematics (and mathematics educators). In general, it seems that the ethics in mathematics (education) literature written by (former) mathematicians is more pragmatic than the writings of education scholars, potentially because the former have to communicate directly with mathematicians. Another example is Buell and Shulman (2021), who write in the abstract of their book on Mathematics for Social Justice, "What makes this book unique and timely is that the most previous curricula linking math and social justice have been treated from a humanist perspective. This book is written by mathematicians, for mathematics students."

For many of the reviewed authors, pragmatism rarely plays an explicit role. However, it can be found when they have personally experienced the tension between (pure) mathematics, the curriculum, students, and teachers themselves and thus had to find practical solutions to particular situations or develop teaching practices that allowed them to navigate these in the long-term. Such pragmatism may also be displayed as a search for ethical guidance in philosophies close to home or one's ancestry (e.g., Abtahi, 2022a, 2022b). For authors working from a more theoretical perspective, pragmatism only seems to play a minor role. This seems particularly true for those whose use of ethics and philosophy is an extension of their existing work on critical mathematics pedagogy.

---

[13] This position of privilege is mentioned here because it has led to conflicts at the Ethics in Mathematics 2 Conference (further explained in Müller (2024)). The question of privilege seems to be an open topic in the discourse on ethics in mathematics itself; particularly, if the turn from sociopolitics to philosophical ethics constitutes a privileged position itself. The question is in so far worthwhile exploring as a critique of privilege may be necessary in a society that aims to overcome injustices (cf., Rieger-Ladich (2023)).

The lesson here is that when it comes to teaching mathematics, a certain pragmatism can be beneficial to navigate the daily "chaos of the classroom." However, too much simplification can be dangerous when it disconnects you from the necessary complexity and ambiguity.[14]

## 4.4. Other positions

We have now discussed three philosophical groundings of ethics in mathematics education. We will now (briefly) consider a final grounding. Much of the literature presented so far critiques Western styles of doing mathematics, and often implicitly, sometimes but not rarely explicitly, of the mainstream culture of Christianity. However, in the spirit of the lessons just derived, these must be considered, even when Christian approaches to ethics in mathematics only represent a minority position in the wider discourse. Let us not forget that the initial scholarship of Freire and others was often supported by Christian institutions (for a historical overview of the developments, see Gottesman, 2016).

---

[14] The tension between complexity and pragmatism/simplification in pedagogy is also explored by Brumfit (1993); the term "chaos of the classroom" is part of the ERIC resume of the Brumfit's article: https://files.eric.ed.gov/fulltext/ED371573.pdf.

Examples include Hansen-Smith's (2003) attempt to bring a Wholeness (i.e., God-centered) ethic into geometry using his (1995, 1999) theory of folding circles. Other works, e.g., Poythress's (1976) biblical view of mathematics is mostly discussed in science and religion studies (e.g., by Howell, 2014) or by scholars not following the modes of critical pedagogy (e.g., Müller, 2018a), and thus often outside of the journals on mathematics education. The discussions appear to struggle with a two-fold hurdle: can and should faith have a bearing on mathematics and arguments against the neutrality of mathematics.

Others explore the ethical dimension of aesthetics in mathematics and the importance of hedonism in mathematics teaching (Lake, 2015, 2019). Both play an important role in what is considered good mathematics among mathematicians. For example, while the 21 characteristics of good mathematics listed by Tao (2007) are comparatively diverse, important aesthetic and hedonic properties still play an important role, including taste, beauty, elegance, depth, and insight. While critical and explicit ethical engagement, as depicted by the scholars discussed in this article, is largely missing in Tao's discussion of these characteristics (Müller, 2024), such discussions may still connect the ethical turn in mathematics education and more conservative mathematics scholars. For example, Peterson (2024) analyzes Serres's (a trained mathematician) perspective on ecology, ethics, and aesthetics, successfully bridging several concerns.

Finally, there are many positions not mentioned here, including ethics of care approaches (e.g., A. Powell & Seed, 2010) and communitarian ethics (e.g., Radford, 2021), which have or can be taken on by individual mathematics education scholars.

Similarly, Abtahi (2021) studies the role of ethics in Vygotsky's zone of proximal development (ZPD), analyzing the power and trust that comes from being the authority responsible for guiding students to more mathematical knowledge. Using aspects from critical mathematics pedagogy and Rorty's (2009) communal view of ethics as sensitivity, she argues for more ethical responsibility by being sensitive towards students. For further reading, we refer the reader to Dubbs (2020).

## 5.  A (brief) critique of the ethical turn

Until now, we have mainly discussed what various authors have gained from considering the ethical dimension of mathematics teaching. Before we move on to some concluding comments, including concrete suggestions for training future teachers, we need to look at some limitations of the existing discourse: empirical data, folk narratives, and the connection to mathematics.

### 5.1.  Empirical scholarship and philosophy

The existing discourse can sometimes build on insufficient empirical studies (as also discussed in Müller, 2024) and argue based on personal experiences and observations. For example, the level-based system proposed by Chiodo and Bursill-Hall (2018) largely comes from observations that the authors made by being in mathematics departments for long enough. The wider literature also *occasionally* writes about mathematicians and mathematics in what could be called folk narratives[15]. Such folk narratives trace back to assumptions or statements made in the past, often with no clear origin, that scholars have taken on over the years. Rarely do these statements get re-evaluated in the discourse. For example, only recently larger empirical studies have been completed about the assumptions and knowledge of ethics by mathematicians (Tractenberg et al., 2024a; Tractenberg et al., 2024b), showing that there might indeed be a desire among many mathematicians to do better, but not necessarily always the knowledge and understanding how to do so. It could be beneficial to examine which narratives are (still) empirically observable and which have turned into folk narratives that have developed over the years and have stuck with the debate for convenience.

---

[15] Here we take inspiration from "folk theorems" in mathematics, as introduced by Rittberg (2023a). Folk theorems are presented as results that are commonly assumed to be true, i.e., as general knowledge, but where either the proof can no longer be found or was never published.

Skovsmose (2020) identifies three competing metanarratives: 1) mathematics as a foundation of civilization, 2) an ally of oppressive forces, and 3) mathematics as an area of (critical) potential for more social justice. All three metanarratives have competing notions of truths, making conversation between them difficult. As Skovsmose points out, these metanarratives can be associated with different political corners and thus have direct impacts on politics and the education system. Our reading of the ethical turn is one that positively extends the sociopolitical turn (narrative 3: mathematics as a critical potential), but we also see a disconnect between what has been written with and without proper empirical grounding. Reducing mathematics education to political narratives or philosophical ethics cannot be the point.

In general, there is an empirical gap when it comes to ethics in mathematics and its education, not necessarily in terms of action ethics in the classroom or sociopolitical questions, but in how mathematicians, teachers, and students perceive, navigate, and understand the relationship between ethics and mathematics. So far, norms in the mathematics classroom have mainly been explored with a socio-mathematical or social focus (Meyer & Schwarzkopf, 2025). The gap, which was also pointed out by Boylan (2016), has meant that discussions found much of their grounding in philosophy rather than data (Ozmantar et al., 2024, p. 2508). Recent empirical studies discuss ethical tensions and ethical dilemmas linked to values, mathematics, and diversity in the classroom (cf. Ozmantar et al., 2024; Roos & Bagger, 2024), with teachers readily using rationalist and non-rationalist approaches to navigate these (Ozmantar et al., 2025). Concrete observable ethical dilemmas include the wide range of mathematical understanding (Alderton & Gifford, 2018), lesson planning and student spontaneity (Foster, 2015), and the exclusion of students' contributions (Brodie, 2010).

However, ethical relevance and awareness also exist outside of ethical dilemmas, and to our knowledge, such situations are only just beginning to be empirically analyzed (e.g., Hauge et al., 2024; Sumpter et al., 2024). Another problem is that mathematical textbooks are rarely analyzed from an ethical rather than sociopolitical perspective (e.g., Dowling, 2024). Few develop ethical theories based on such observations (e.g., Roth, 2024), and the situation is made more difficult as fewer teachers seem to be writing about ethics than sociopolitical issues, so despite all good desires, the practicality of the ethical turn, is yet another question to be solved (cf. Pais, 2024).

This gap may not necessarily be restricted to ethics in mathematics education or its empirical dimension. For example, Gottesman's (2016) historical study of the critical turn in education showed that later authors sometimes tended to take on board broad narratives and philosophical assumptions (e.g., (neo-)Marxist assumptions about progress and conflict) – quite often unknowingly and particularly from influential founding figures (e.g., Paulo Freire): "Too often it feels like we do not read ideas closely enough. Much of the problem, I believe, is that we rely too much on secondary sources." (ibid., 2016, p. 169)

So far, the results of the few existing empirical studies support the largely philosophical discourse about the necessity to navigate ethics in everyday classroom interactions, but much remains to be empirically validated.

## 5.2. Connected to and disconnected from mathematics

Large parts of the literature on ethics in mathematics education build on critical mathematics pedagogy and thus inherently take on many of the normative assumptions of critical theory, especially the emancipatory goal and the critique of existing norms and value systems of society. Such a commitment to a better world shapes many of the reviewed authors' scientific and educational approaches. This can clash with the science of mathematics, which may understand itself as neutral, pure, universal, rational, and ultimately true.

Such dynamics are often found in discourses building on Postmodern, poststructuralist, and postcolonial research, where "scholars have provincialized or localized [the] universality claim [of knowledge]. For them, the claim to universality and foundation is only a question of epistemic power. The universal and the foundational knowledge have been questioned as a God trick of scholars […] [and] moreover, the relativism has been redefined not only as relative to time and place (historical reason) but also to irreducible and incommensurable differences of societies and cultures (epistemological reason)" (Canceran, 2023, p. 8; building on the works of Nagel (1989) and Chakrabarty (2000)). However, while this is a common phenomenon in the interaction between historians, philosophers, and educators of science with scientists, it is particularly hard-hitting in the discourse on ethics in mathematics.

While, at times, this means that the discourse shows little common ground between the moral issues that mathematicians may see and those that education scholars observe, some try to change this by trying to act as a bridge between the separate disciplines and perspectives (e.g., Chiodo & Müller, 2024a; Ernest, 2021c; Müller, 2024), among others by reframing an initially social justice-oriented motivation for doing mathematics into questions of (philosophical) ethics (cf., Piercey, 2023, p. 66), and by trying to bridge different parts of philosophy (Paton & Sinclair, 2024) or by bringing school mathematics and research mathematics closer together (Sarikaya, 2024). Reframing issues from social justice to ethics, and thus philosophy, may make it easier for mathematicians to connect with the discourse.

However, there is little reflection on the meta-level of how such a switch can support a deeper conversation between mathematicians and educators. At present, the sociopolitical and ethical discourses present some convergent dynamics, as presented in (Tractenberg et al., 2024a; Tractenberg et al., 2024b), and it also shows clear divergences (cf. Crowell, 2022; Rittberg, 2024). The ethical turn reflects on the sociopolitical questions, but to be truly successful, it may need to reflect more on itself, too.

## 6. <u>Training (future) teachers</u>

Ethical issues in the classroom need not always be as obvious as the exercise[16] that a German textbook posed during World War 2:

"Problem 200. According to statements of the Draeger Works in Luebeck in the gassing of a city only 50% of the evaporated poison gas is effective. The atmosphere must be

---

[16] Shulman (2002) also discussed an economics book from 1974 that contained a question on slave whipping, and a separate optimisation problem for finding the optimal position of an offshore oil platform.

poisoned up to a height of 20 metres in a concentration of 45 mg. per cubic metre. How much phosgene is needed to poison a city of 50,000 inhabitants who live in an area of four square kilometres? How much phosgene would the population inhale with the air they breathe in ten minutes without protection against gas, if one person uses 30 litres of breathing air per minute? Compare this quantity with the quantity of the poison gas used." (Shulman, 2002, pp. 123–124; referencing Cohen, 1988, p. 244)

In today's mathematics classrooms, they are often much more subtle. Our exploration has revealed a multifaceted literature on ethics in mathematics education that largely mirrors the complexity that can be present in (urban) classrooms (cf. Rivera, 1998). At times, the discourse is marked by a tension between "doing good" and "preventing harm." At other times, there is an unforgiving tension between the "should" and the "can" (especially when public policy stipulates a different form of teaching than is advocated by the surveyed literature; a conflict between reality and literature that is present in many areas of mathematics education, see also (Fiona, 2005)). On the one hand, scholars aim to leverage ethics in mathematics to foster democratic participation, empower students, and even advance their emotional development. On the other, a significant strand also focuses on reducing and mitigating potential negative consequences of mathematics, ranging between Ernest's (2016a) "collateral damage" of learning mathematics and O'Neil's (2016) "weapons of math destruction."

Philosophical assumptions and theories about the nature of mathematics can influence didactical decisions(Ernest, 1989; Maaß & Götz, 2022), and this can lead to a perspective on ethics in mathematics education where the subject (i.e., mathematics) is supposedly neutral and the education is supposedly supposed to be neutral, too. However, "neutrality is not a solution" (Wohnig & Zorn, 2022) in civic education, and agreements such as the Beutelsbach consensus, can, and maybe should, be read as a controversy demand rather than a controversy ban (Frick, 2022) – a challenge of increasing importance when ethics, sustainability and mathematics education meet.[17]

Various potential avenues for teaching ethics in mathematics have been explored in the literature and practice, including embedded ethics (Chiodo et al., 2023; Müller & Chiodo, 2023; Shah, 2022), mini-seminars (Miller, 2022), regular examined seminars (Cordes, 2023), regular non-examinable seminars (CUEiMS, 2016), courses (Franklin, 2005; Skufca, 2021), card games (Walk & Tractenberg, 2024) for use in professional development workshops, and the outsourcing of ethics training (Alayont, 2022, p. 167). However, specific ethics courses aimed at mathematics students remain comparatively rare (as also discussed in Ernest, 2024c) and might not be sufficient to normalize thinking about ethics in mathematics and its education (Chiodo & Bursill-Hall, 2019).

When one widens the context to sociopolitical and social justice concerns, more has been tried and tested, but it is still not the mainstream in mathematics. While the situation looks better in mathematics education, especially regarding sociopolitical issues, explicit ethics in mathematics education courses remain rare, especially courses beyond research ethics and education ethics. However, how much explicit training on ethics in mathematics future teachers need is also still an open question, but the motivation seems to be there.[18] Regarding ethics, some recent studies from adjacent

---

[17] Also explored by Wilhelm (2024) by building on Andelfinger's thoughts on gentle mathematics education.

[18] A study performed by Karatas and Yilmaz (2021) on future Turkish mathematics teachers remarked "that the prospective teachers showed an interest in ethical values, such as justice, being a good role model and honesty. They describe the ethical teacher as being 'of holy personality, a guide and professional', and the unethical teacher as 'useless and harmful'." (ibid., p. 377)

fields[19] suggest that it may be less than expected. Some older research also suggests a nuanced picture (Colnerud, 1997).

The lessons we draw from the ethical turn in mathematics education can be summarized as expanding teachers' professional identity to include self-critical personal development and an inward journey[20] towards seeing the "richness of being"[21] in the mathematics classroom. As the ethical turn strongly focuses on the interrelations of individual people in and outside the classroom, the ethical awareness of the teacher becomes central to turning theoretical scholarship into praxis. Mathematics is seen as much more than just knowledge; it requires solidarity with people, animals, and nature (as also explored by Abtahi, 2024). As such, we argue for the additional training of future teachers in the following seven areas:

1. _Knowing the self_ and being open to self-criticism is foundational to ethical relationships in the classroom.

2. _Accepting the Other_ by recognizing each student's unique worth and needs, fostering inclusive classroom environments, and taking on personal responsibility for engaging in ethical teacher-student interactions.

3. _Studying some philosophy and history of mathematics_ empowers teachers to question taken-for-granted assumptions and present mathematics as a human practice.

---

[19] For a study of Turkish technology teachers, we refer to Baysan and Cetin (2021), who found that the vast majority of teachers had a good ethical basis.

[20] An example of this is given by House (2024).

[21] As also analysed by Maheux (2024) in his discussion of Feyerabend.

4. *Appreciating the societal role of mathematics* fosters a wider sense of awareness and responsibility for how mathematics is portrayed and used in and outside the classroom.

5. *Being flexible and pragmatic* is vital to adapting abstract ethical principles to the dynamic realities of teaching.

6. *Looking beyond one's culture* enables connecting with today's diverse classrooms more effectively.

7. *Developing bravery and preparedness* to develop courage and willingness to challenge the mainstream ways mathematics is practiced, taught, and examined.

We note that these seven suggestions are closely linked to the typical categories of ethical conflict in the classroom, i.e., interpersonal norms, internal professional norms, institutional norms, social conformity norms, and self-protecting norms (Colnerud, 1997, p. 630). Fluency in ethics can influence teachers' identities (Granjo et al., 2021), and the lessons coming out of the discourse on ethics in mathematics education take aim at this. Nevertheless, for a more comprehensive picture, they will likely need to be accentuated by the wider lessons from the sociopolitical turn and professional ethics in education (e.g., the focus on specific aspects such as equity, honesty, and (intellectual) freedom).

In short, a careful reading of the ethical turn in mathematics education suggests that the ethical training of mathematics teachers is both a matter of professional ethics and (individual) moral education. Professional ethics likely falls on empty ground without the latter. The literature puts humans at the center of mathematics while trying to advocate for a less anthropocentric view and usage of mathematics, while at times trying to move beyond[22] a postmodern analytic perspective. These lessons thus put the teacher as a person into the center of attention rather than a specific vision for mathematics.

## 7. <u>Conclusion</u>

As explored throughout this article, the emerging ethical turn is not a "turning away" from sociopolitical questions but a nuanced deepening and extension. While the latter provides the frameworks for understanding the social, cultural, and political nature of mathematics, the ethical turn brings a sharpened focus onto the individual – the learner, the teacher, the mathematician – and their inherent responsibilities within the practice and teaching of mathematics. This shift from sociopolitics to ethics involves an explicit engagement with (philosophical) ethics, a strengthened challenge of the neutrality, purity, and universality of mathematics, attempts to bridge different mathematical communities, and a balance between intellectual and activist scholarship.

Over the years, the discourse on ethics in mathematics education has shifted from a limited early engagement with ethics before 2000 towards a deeper engagement with ethics in recent years. The early reception of Levinas was a notable turning point. Quickly, social justice concerns were integrated with ethics in mathematics, and the first multidimensional ethical frameworks emerged. Today, the discourse focuses on combining philosophical insight with practical implementations, and we are seeing the beginning of empirical studies with an explicit focus on ethics in the mathematics classroom.

---

[22] Given the deep roots of the ethical turn in postmodern ideas and scholarship, how exactly ethics in mathematics education can move beyond it, is hard to ascertain. What this means for mathematics education more general has also been explored by Skovsmose (2012).

This growing ethical awareness among mathematics educators carries important insights for the training of future mathematics teachers. It urges us to move beyond technical and pedagogical training towards understanding mathematics in a larger sociopolitical and ethical context to cultivate a moral agency that is grounded in theory, well-reflected, and sufficiently strong for the classroom. Navigating the many dimensions of mathematics and the complex tension between using mathematics to "do good" and the need to "prevent harm" requires us to look beyond professional ethics towards an ethical Self of the teacher, including fostering self-reflection, responsibility towards the Other, historical and philosophical understanding, flexibility and pragmatism, cultural sensitivity, and the courage to address ethical dilemmas.